# Number of Partitions of an $n$-kilogram Stone into Minimum Number of Weights to Weigh All Integral Weights from 1 to $n$ kg(s) on a Two-pan Balance

**Md Towhidul Islam**

Assistant Professor

Department of Marketing, Comilla University

Courtbari, Comilla

Bangladesh

towhid@cou.ac.bd

towhid81@gmail.com

**Md Shahidul Islam**

Deputy Director

Bangladesh Railway, Bangladesh Civil Service

Dhaka, Bangladesh

shahid1514@gmail.com

*Abstract*

*We find out the number of different partitions of an $n$-kilogram stone into the minimum number of parts so that all integral weights from 1 to $n$ kilograms can be weighed in one weighing using the parts of any of the partitions on a two-pan balance. In comparison to the traditional partitions, these partitions have advantage where there is a constraint on total weight of a set and the number of parts in the partition. They may have uses in determining the optimal size and number of weights and denominations of notes and coins.*

## Introduction

A seller has an $n$ kilogram stone which he wants to break into the minimum possible number of weights using which on a two-pan balance he can sell in whole kilograms up to $n$ kilogram(s) of goods in one weighing. As in tradition, he can place weights on both the pans but goods on only one pan. We call such a partition a 'feasible' partition. Our intension is to find out the number of all feasible partitions of $n$.

Suppose the minimum possible $m$ weights are $w_1 \leq w_2 \leq \cdots \leq w_{m-1} \leq w_m$ and $n = w_1 + w_2 + \cdots + w_{m-1} + w_m$. If we also suppose $u_i = [-1, 0, 1]$, from the description of the problem it can be stated that $u_1 w_1 \mp u_2 w_2 \mp \cdots \mp u_{m-1} w_{m-1} \mp u_m w_m$ must make all the positive integers from 1 to $n$. Then we will find out $t(n)$, the number of all such partitions of the ordered integral set of weights $w_1 \leq w_2 \leq \cdots \leq w_{m-1} \leq w_m$ for every positive integer $n$. To that end we at first define some terms used in the discussion. Then we prove a number of theorems which ultimately lead to our main finding, the recursive functions for $t(n)$ in two spans. O'Shea [1] introduced the concept of M-partitions by partitioning a weight into as few parts as possible so as to be able to weigh any integral weight less than $m$ weighed on a one scale pan. He maintained the subpartition property of MacMahon's [2] perfect partitions but dropped the uniqueness property and also added a new property; the minimality of number of parts in the partitions. We examine the partitioning situation for a two-pan balance maintaining his minimality of parts.

## Definitions

A *feasible set/partition* of $n$ is an ordered partition of minimum possible $m$ parts $w_1 \leq w_2 \leq \cdots \leq w_{m-1} \leq w_m$ made from $n$ such that all integral values from 1 to $n$ can be weighed in one weighing using the parts on a traditional two-pan balance.

$R_i = w_1 + w_2 + \cdots + w_i$. From this, it is clear that $R_1 = w_1$, $R_i = R_{i-1} + w_i$ and $R_m = w_1 + w_2 + \cdots + w_m = n$. We assume $R_0 = 0$.

$t(n)$ is the number of all feasible partitions of $n$. We assume $t(0) = 1$.

## Theorems of Feasibility

**Theorem 1.** *For any feasible set, the lightest weight $w_1$ equals 1 kg, i.e. $w_1 = 1$*

**Proof.** Suppose $w_1 \neq 1$. So, $w_i \geq 2$ for all $i$. However, using such a partition, $n - 1$ kg cannot be weighed. Placing all $w_i$ pieces on one of the pans we see, $w_1 + w_2 + \cdots + w_{m-1} + w_m = n$. Now, if we take out the smallest piece $w_1$ from there, we see $w_2 + w_3 \ldots + w_{m-1} + w_m \leq n - 2$. So, such a partition with $w_1 \neq 1$ can never weigh $n - 1$. Therefore, we conclude the lightest weight $w_1 = 1$ kg.

**Theorem 2.** *For any feasible partition, $w_i \leq 2R_{i-1} + 1$.*

**Proof.** Two things are clear-
i) The highest value possible to be weighed from $u_1 w_1 \mp u_2 w_2 \mp \cdots \mp u_{i-1} w_{i-1}$ is $w_1 + w_2 + \cdots + w_{i-1} = R_{i-1}$.
ii) So, the lowest possible value made from $w_i - (u_1 w_1 \mp u_2 w_2 \mp \cdots \mp u_{i-1} w_{i-1})$ is $w_i - R_{i-1}$.

Therefore, for a feasible set, there should not be any integer $n$ in the range $R_{i-1} < n < w_i - R_{i-1}$.

That is, $R_{i-1} + 1 \geq w_i - R_{i-1}$.

$\Rightarrow w_i \leq 2R_{i-1} + 1$.

***Corollaries of Theorem 2.*** *Interrelationships among* $w_i$, $R_i$ *and* $R_{i-1}$

From theorem 2, we get $\quad R_{i-1} \geq \frac{w_i - 1}{2}$ .... Corollary (1)

$\Rightarrow R_i - w_i \geq \frac{w_i - 1}{2}$ (From definition, $R_{i-1} = R_i - w_i$)

$\Rightarrow R_i \geq \frac{3w_i - 1}{2}$ .... Corollary (2)

$\Rightarrow w_i \leq \frac{2R_i + 1}{3}$ .... Corollary (3)

$\Rightarrow R_i - R_{i-1} \leq \frac{2R_i + 1}{3}$ (From definition, $R_i - R_{i-1} = w_i$)

$\Rightarrow R_i \leq 3R_{i-1} + 1$ .... Corollary (4)

$\Rightarrow R_{i-1} \geq \frac{R_i - 1}{3}$ .... Corollary (5)

**Theorem 3.** *The highest ever possible value of part* $w_i$ *is* $3^{i-1}$ *and the highest ever possible value of* $R_i$ *feasibly partitionable in i parts is* $\frac{3^i - 1}{2}$.

**Proof.** From Theorem 2, we know $w_i \leq 2R_{i-1} + 1$ for $2 \leq i \leq m$; from definition we know $R_1 = w_1$ and from Theorem 1 we know, $w_1 = 1$.

So, $w_2 \leq 2R_1 + 1$

$\Rightarrow w_2 \leq 3$

So, the highest possible value of $R_2$ is $w_1 + w_2 = 1 + 3 = 4$

Going on with $w_i \leq 2R_{i-1} + 1$, we see the highest possible values of weights $w_1, w_2, w_3 ..., w_m$ are $3^0, 3^1, 3^2, ..., 3^{m-1}$. Clearly, $w_i = 3^{i-1}$.

And the highest possible $R_i$ feasibly partitionable in $i$ pieces is $R_i = 3^0 + 3^1 + \cdots + 3^{i-2} + 3^{i-1} = \frac{3^i - 1}{2}$.

**Theorem 4.** *At least m weights are needed for* $\frac{3^{m-1} + 1}{2} \leq n \leq \frac{3^m - 1}{2}$ *where* $m = \lceil log_3(2n) \rceil$.

**Proof.** From Theorem 3, it is clear that he highest possible value of $R_m$ feasibly partitionable in $m$ pieces is $R_m = \frac{3^m - 1}{2}$ and the highest possible value of $R_{m-1}$ feasibly partitionable in $m - 1$ pieces is $R_{m-1} = \frac{3^{m-1}-1}{2}$. So, the next integer $\frac{3^{m-1}-1}{2} + 1 = \frac{3^{m-1}+1}{2}$ is the lowest value of $n$ feasibly partitionable in $m$ parts.

So, we have proved that at least $m$ weights are needed for $\frac{3^{m-1}+1}{2} \leq n \leq \frac{3^m-1}{2}$.

$$\Rightarrow 3^{m-1} + 1 \leq 2n \leq 3^m - 1$$
$$\Rightarrow 3^{m-1} < 2n < 3^m$$
$$\Rightarrow m - 1 < log_3(2n) < m$$

However, $m$ is never a fraction. The ceiling function $m = \lceil log_3(2n) \rceil$ always gives the $m$ we know to be the correct number of parts for the range $\frac{3^{m-1}+1}{2} \leq n \leq \frac{3^m-1}{2}$.

**Theorem 5.** For $\frac{3^{m-1}+1}{2} \leq n \leq \frac{3^{m-1}+1}{2} + 3^{m-2}$, the range of $R_{m-1}$ is $\lceil \frac{n-1}{3} \rceil \leq R_{m-1} \leq \lfloor \frac{2n+3^{m-2}-1}{4} \rfloor$ and for $\frac{3^{m-1}+1}{2} + 3^{m-2} + 1 \leq n \leq \frac{3^m-1}{2}$, the range of $R_{m-1}$ is $\lceil \frac{n-1}{3} \rceil \leq R_{m-1} \leq \frac{3^{m-1}-1}{2}$.

**Proof.**

**The smallest feasible $R_{m-1}$**

Putting $i = m$ in Corollary 5 of Theorem 2 we get, $R_{m-1} \geq \frac{R_m-1}{3}$

As $R_m = n$, the smallest feasible $R_{m-1}$ is $\lceil \frac{n-1}{3} \rceil$ for all $n$.

For example, for $n = 16$, the smallest feasible $R_{m-1} = \lceil \frac{n-1}{3} \rceil = 5$ and for $n = 26$, the smallest feasible $R_{m-1} = \lceil \frac{n-1}{3} \rceil = 9$.

**The largest feasible $R_{m-1}$**

With a little modification of Theorem 3, it is clear that for any $n$ the $m$-part partition $3^0, 3^1, ..., 3^{m-3}, \lfloor \frac{n-R_{m-2}}{2} \rfloor, \lceil \frac{n-R_{m-2}}{2} \rceil$ would ensure the highest possible $R_{m-1}$ given the condition that $\lceil \frac{n-R_{m-2}}{2} \rceil \leq 3^{m-2}$.

Based on whether $\lceil \frac{n-R_{m-2}}{2} \rceil \leq 3^{m-2}$ or not, we can split the range $\frac{3^{m-1}+1}{2} \leq n \leq \frac{3^m-1}{2}$ found in Theorem 4 into two mutually exclusive and collectively exhaustive spans: a) $\frac{3^{m-1}+1}{2} \leq n \leq \frac{3^{m-1}+1}{2} + 3^{m-2}$ and b) $\frac{3^{m-1}+1}{2} + 3^{m-2} + 1 \leq n \leq \frac{3^m-1}{2}$.

a) For $\frac{3^{m-1}+1}{2} \leq n \leq \frac{3^{m-1}+1}{2} + 3^{m-2}$, we see $\left\lceil \frac{n-R_{m-2}}{2} \right\rceil \leq 3^{m-2}$. So, the largest feasible $R_{m-1}$ is $3^0 + 3^1 + \cdots + 3^{m-3} + \left\lceil \frac{n-R_{m-2}}{2} \right\rceil = \left\lceil \frac{2n+3^{m-2}-1}{4} \right\rceil$.

For example, if $n = 16$, the largest feasible $R_{m-1}$ is $\left\lceil \frac{2n+3^{m-2}-1}{4} \right\rceil = 10$.

b) However, for $\frac{3^{m-1}+1}{2} + 3^{m-2} + 1 \leq n \leq \frac{3^m-1}{2}$, the condition is not met. From Theorem 3 we know $w_{m-1}$ should never exceed $3^{m-2}$. So the $m$-part partition should be $3^0, 3^1, \ldots, 3^{m-3}, 3^{m-2}, n - \frac{3^{m-2}-1}{2}$. And the largest feasible $R_{m-1}$ is $3^0 + 3^1 + \cdots + 3^{m-3} + 3^{m-2} = \frac{3^{m-1}-1}{2}$.

For example, if $n = 26$, the largest feasible $R_{m-1}$ is $\frac{3^{m-1}-1}{2} = 13$.

So, we have found the range of $R_{m-1}$ for the two segments as described in the statement of this theorem.

## The Main Result: The Recursive Functions for $t(n)$

**Theorem 6.**

$$t\left(\frac{3^{m-1}+1}{2} \leq n \leq \frac{3^{m-1}+1}{2} + 3^{m-2}\right)$$
$$= \sum_{R_{m-1}=\left\lceil \frac{n-1}{3} \right\rceil}^{\left\lceil \frac{2n+3^{m-2}-1}{4} \right\rceil} t(R_{m-1}) - \sum_{R_{m-1}=\left\lceil \frac{3n+2}{5} \right\rceil}^{\left\lceil \frac{2n+3^{m-2}-1}{4} \right\rceil} \sum_{R_{m-2}=\left\lceil \frac{R_{m-1}-1}{3} \right\rceil}^{2R_{m-1}-n-1} t(R_{m-2})$$

and

$$t\left(\frac{3^{m-1}+1}{2} + 3^{m-2} + 1 \leq n \leq \frac{3^m-1}{2}\right) = \sum_{R_{m-1}=\left\lceil \frac{n-1}{3} \right\rceil}^{\frac{3^{m-1}-1}{2}} t(R_{m-1})$$

**Proof.** In order to determine the number of feasible partitions of $n$ we will derive two different recursive functions for the two spans from Theorem 5.

a) *Determining $t(n)$ for the range* $\frac{3^{m-1}+1}{2} \leq n \leq \frac{3^{m-1}+1}{2} + 3^{m-2}$

To count the total number of feasible partitions of $n$, we have to count all feasible partitions of all $R_{m-1}$ possible to be broken from $n$ because adding an additional last part $w_m$ to these $m-1$ part feasible partitions of $R_{m-1}$ will turn them into $m$ part feasible partitions of $n$; the condition is that no $w_{m-1}$ possible to be made from each of these $R_{m-1}$ values exceeds the corresponding $w_m$.

From Theorem 5 we know, for $\frac{3^{m-1}+1}{2} \leq n \leq \frac{3^{m-1}+1}{2} + 3^{m-2}$ the range of $R_{m-1}$ is $\left\lceil \frac{n-1}{3} \right\rceil \leq R_{m-1} \leq \left\lfloor \frac{2n+3^{m-2}-1}{4} \right\rfloor$ and $w_m = \left\lceil \frac{n-R_{m-2}}{2} \right\rceil$. So, the formula at first seems to be

$$t\left(\frac{3^{m-1}+1}{2} \leq n \leq \frac{3^{m-1}+1}{2} + 3^{m-2}\right) = \sum_{R_{m-1}=\left\lceil \frac{n-1}{3} \right\rceil}^{\left\lfloor \frac{2n+3^{m-2}-1}{4} \right\rfloor} t(R_{m-1})$$

However this would count some partitions in duplication as for some values of $n$ in this range, some of the $w_{m-1}$ values feasibly broken from each of these $R_{m-1}$ values are larger than the corresponding $w_m = \left\lceil \frac{n-R_{m-2}}{2} \right\rceil$. If arranged in ascending order, it would be clear that the partitions are counted in duplication for some other $n$. We have to exclude those duplications from the count by finding out such partitions of these $R_{m-1}$ values.

To do that, we will at first find the range of such problematic $R_{m-1}$ values and then we will set the range of incompatible $R_{m-2}$ values for each of these $R_{m-1}$ values.

From corollary 3 of Theorem 2 it is clear that $w_{m-1} \leq \frac{2R_{m-1}+1}{3}$ and from definition, $w_m = n - R_{m-1}$. It is noticeable that $w_{m-1}$ will be greater than $w_m$ if $w_m \leq w_{m-1} - 1$.

$\Rightarrow n - R_{m-1} \leq \frac{2R_{m-1}+1}{3} - 1$.

$\Rightarrow \frac{3n+2}{5} \leq R_{m-1}$.

This lower limit of $R_{m-1}$ taken in consideration along with the range of $R_{m-1}$ set in Theorem 5 redefines the range of such problematic $R_{m-1}$ as $\frac{3n+2}{5} \leq R_{m-1} \leq \left\lfloor \frac{2n+3^{m-2}-1}{4} \right\rfloor$ where some of the $w_{m-1}$ values broken from $R_{m-1}$ are larger than the corresponding $w_m$.

As $w_{m-1}$ must not be larger than $w_m$, we have to exclude those partitions of these $R_{m-1}$ in $m-1$ parts where $w_m + 1 \leq w_{m-1} \leq \frac{2R_{m-1}+1}{3}$. (As from corollary 3 of Theorem 2 we know, $w_{m-1} \leq \frac{2R_{m-1}+1}{3}$.)

That is, we have to exclude partitions with $n - R_{m-1} + 1 \leq w_{m-1} \leq \frac{2R_{m-1}+1}{3}$

Or, $R_{m-1} - n + R_{m-1} - 1 \geq R_{m-1} - w_{m-1} \geq R_{m-1} - \frac{2R_{m-1}+1}{3}$ (deducting the terms from $R_{m-1}$).

Or, $\frac{R_{m-1}-1}{3} \leq R_{m-2} \leq 2R_{m-1} - n - 1$ for each of these problematic $R_{m-1}$.

Finally the formula stands,

$$t\left(\frac{3^{m-1}+1}{2} \leq n \leq \frac{3^{m-1}+1}{2} + 3^{m-2}\right)$$

$$= \sum_{R_{m-1}=\left\lceil\frac{n-1}{3}\right\rceil}^{\left\lfloor\frac{2n+3^{m-2}-1}{4}\right\rfloor} t(R_{m-1}) - \sum_{R_{m-1}=\left\lceil\frac{3n+2}{5}\right\rceil}^{\left\lfloor\frac{2n+3^{m-2}-1}{4}\right\rfloor} \sum_{R_{m-2}=\left\lceil\frac{R_{m-1}-1}{3}\right\rceil}^{2R_{m-1}-n-1} t(R_{m-2})$$

b) *Determining $t(n)$ for the range* $\frac{3^{m-1}+1}{2} + 3^{m-2} + 1 \leq n \leq \frac{3^m-1}{2}$

From Theorem 5 we know, for $\frac{3^{m-1}+1}{2} + 3^{m-2} + 1 \leq n \leq \frac{3^m-1}{2}$ the range of $R_{m-1}$ is $\left\lceil\frac{n-1}{3}\right\rceil \leq R_{m-1} \leq \frac{3^{m-1}-1}{2}$. However, unlike for the range of $n$ in part (a), all possible $w_{m-1}$ values broken from $R_{m-1}$ in this range are less than the corresponding $w_m$. So, to count the total possible number of feasible partitions of $n$ in this range, we have only to count $t(R_{m-1})$ for all possible $R_{m-1}$ values for $n$; no chance of duplication arises. So, the formula stands,

$$t\left(\frac{3^{m-1}+1}{2} + 3^{m-2} + 1 \leq n \leq \frac{3^m-1}{2}\right) = \sum_{R_{m-1}=\left\lceil\frac{n-1}{3}\right\rceil}^{\frac{3^{m-1}-1}{2}} t(R_{m-1})$$

When we find out the $t\left(\frac{3^{m-1}+1}{2} + 3^{m-2} + 1 \leq n \leq \frac{3^m-1}{2}\right)$ values starting from $\frac{3^m-1}{2}$ backwards, we see the terms of sequence A005704 of the OEIS [3] come in triplicates.

## References

[1] E. O'Shea, M-partitions: Optimal Partitions of Weight for One Scale Pan, Discrete Math. 289 (2004), 81–93.

[2] P. A. MacMahon, The Theory of Perfect Partitions and the Compositions of Multipartite Numbers, Messenger of Math. 20, (1891), 103–119.

[3] The Online Encyclopedia of Integer Sequence, at https://oeis.org/A005704.